\documentclass[12pt]{article}
\usepackage{amssymb,amsmath,amscd,epsfig,amsthm}
\usepackage{dsfont}

\newcommand{\Stg}{\mathop{\rm Stab}\nolimits_\G}
\newcommand{\Stga}{\mathop{\rm Stab}\nolimits_\Gamma}
\newcommand{\G}{\mathcal G}
\newcommand{\A}{\mathcal A}
\newcommand{\B}{\mathcal B}
\newcommand{\one}{\mathds 1}
\newcommand{\zero}{\mathds O}
\newcommand{\Sym}{\mathop{\rm Sym}\nolimits}
\newcommand{\Gamman}{\Gamma^{(n)}}

\newcommand{\Tn}{T^{(n)}}

\newtheorem{theorem}{Theorem}[section]
\newtheorem{cor}[theorem]{Corollary}
\newtheorem{prop}[theorem]{Proposition}
\newtheorem{lemma}[theorem]{Lemma}
\newtheorem{defin}{Definition}

\begin{document}
\title{Automata generating free products\\ of groups of order 2}
\author{Dmytro Savchuk\thanks{Supported by NSF grants DMS-0308985 and
DMS-0456185}, Yaroslav Vorobets\\ Texas A\&M University}

\maketitle

\abstract{We construct a family of automata with $n$ states, $n\geq
4$, acting on a rooted binary tree that generate the free products
of cyclic groups of order 2.}

\section*{Introduction}
An automaton group is the group generated by transformations defined
by all states of a finite invertible automaton (see precise
definitions in Section~\ref{sec_prelim}) over a finite alphabet $X$.
The transformations act on the set $X^*$ of finite words over $X$,
which can be regarded as a regular rooted tree. The first mentioning
of automaton groups dates back to early
1960s~\cite{glushkov:ata,horejs:automata}. Until the beginning of
1980s the interest in these groups was somewhat sporadic. It started
to grow rapidly after it was shown that the class of automaton
groups contains counterexamples to the general Burnside problem
\cite{aleshin:burnside,sushch:burnside,grigorch:burnside,gupta_s:burnside}.
Later, Grigorchuk in~\cite{grigorch:degrees} solved Milnor problem
on the existence of groups of intermediate growth by producing
series of groups generated by automata. A very abridged list of
automaton groups that have extraordinary properties includes:
amenable, but not elementary amenable group
(see~\cite{grigorch:example}) answering Day
problem~\cite{day:amenable}; amenable, but not subexponentially
amenable Basilica
group~\cite{grigorch_z:basilica,bartholdi_v:amenab}; lamplighter
group that gave rise to the negative solution of the strong Atiyah
conjecture on $L_2$-Betti numbers~\cite{grigorch_lsz:atiyah}.

All transformations defined by states of finite invertible automata
over a fixed alphabet form a group of automatic transformations over
this alphabet. The structure of this large group is yet to be
understood. An interesting question is the embedding of known groups
into this group. For example, Brunner and Sidki proved
in~\cite{brunner_s:glnz} that $GL_n(\mathbb Z)$ can be generated by
finite automata over the alphabet with $2^n$ letters. In this paper
we address this question in regard to the free products of groups of
order 2 (we will often denote the group of order 2 by $C_2$). The
first embedding of such free products into the group of automatic
transformations over the 2-letter alphabet was constructed by
Olijnyk~\cite{olijnyk:C2}. We also mention results of C.~Gupta,
N.~Gupta and
A.~Olijnyk~\cite{olijnyk:free_products,gupta_go:free_products} who
embedded the free product of any finite family of finite groups into
a group of automatic transformations over a suitable alphabet.

The above constructions lack the important property of
self-similarity~\cite{nekrash:self-similar}. In other words, the
group is not generated by all states of a single automaton. The
first self-similar example was provided by a 3-state automaton
$\B_3$ over 2-letter alphabet whose Moore diagram is depicted in
Figure~\ref{fig_bel}. This automaton was studied during the summer
school in Automata groups held in 2004 at the Autonomous University
of Barcelona in Bellaterra. Since then, it is known as the
Bellaterra automaton. It was proved by Muntyan (see the proof
in~\cite{bondarenko_gkmnss:full_clas32}
or~\cite{nekrash:self-similar}) that $\B_3$ generates the group
isomorphic to the free product of 3 copies of groups of order 2.

Many papers on free groups and free products generated by automata
share the same idea of dual automaton. For an automaton $\A$ the
dual automaton $\hat \A$ is obtained from $\A$ by interchanging the
states and the alphabet, and swapping the transition and output
functions. For precise definition see Section~\ref{sec_prelim}. It
turns out that the ``freeness'' properties of the group generated by
$\A$ are related to certain transitivity conditions of the action of
the group generated by $\hat\A$.

The Bellaterra automaton belongs to the class of bireversible
automata~\cite{nekrash:self-similar,gl_mo:compl}, which seems to be
a natural source for automata generating free groups and free
products. An invertible automaton is called bireversible if its dual
and the dual to its inverse are also invertible. It is worth
mentioning that the Bellaterra automaton was discovered while
classifying all bireversible 3-state automata over 2-letter
alphabet.

The transformations $a,b,c$ defined by states of the Bellaterra
automaton act on the set $\{0,1\}^*$ of finite words over the
alphabet $X=\{0,1\}$. They are uniquely determined by so-called
wreath recursion
\begin{equation*}
\begin{array}{lcl}
a&=&(c,b),\\
b&=&(b,c),\\
c&=&(a,a)\sigma,\\
\end{array}
\end{equation*}
where $t=(t_0,t_1)$ means that $t(0w)=0t_0(w)$ and $t(1w)=1t_1(w)$
for any word $w\in X^*$ while $t=(t_0,t_1)\sigma$ means
$t(0w)=1t_0(w)$ and $t(1w)=0t_1(w)$ for any $w\in X^*$.

The Bellaterra automaton gives rise to a family of bireversible
automata in which all states define involutive transformations. The
construction is very simple. Namely, we modify the automaton $\B_3$
by inserting new states on the path from $c$ to $a$. More precisely,
each automaton in the family is defined by wreath recursion
\begin{equation}
\label{eqn_gen_family}
\begin{array}{lcl}
a&=&(c,b),\\
b&=&(b,c),\\
c&=&(q_{1},q_{1})\sigma_{0},\\
q_{i}&=&(q_{i+1},q_{i+1})\sigma_{i},\ i=1,\ldots,n-4,\\
q_{n-3}&=&(a,a)\sigma_{n-3},\\
\end{array}
\end{equation}
where $\sigma_i\in\Sym(\{0,1\})$ is chosen arbitrarily.

Conjecturally, each  automaton in the family for which at least one
of the $\sigma_i$ is nontrivial, generates the free product of
groups of order 2. The first result supporting this conjecture was
obtained by M.~Vorobets and Y.~Vorobets~\cite{vorobets:series_free}.
It was shown that if the number of states is odd and $\sigma_i=(12)$
for all $i$, then the conjecture holds. In the subsequent paper by
the same authors and B.~Steinberg~\cite{steinberg_vv:series_free}
the conjecture was proved for the automata with even number of
states and additional condition that the number of nontrivial
$\sigma_i$ is odd.

In this paper we prove that any $n$-state automaton from the
family~\eqref{eqn_gen_family} with $n\geq 4$ satisfying
$\sigma_0=(12)$ and $\sigma_{n-3}=(12)$ generates the free product
of groups of order 2. This result covers the series constructed
in~\cite{vorobets:series_free} except one, but the most important
case $n=3$, and partially overlaps with a family constructed
in~\cite{steinberg_vv:series_free}. More precisely, our main results
are as follows.

Let $\B_4$ be the automaton defined by the wreath recursion
\begin{equation*}
\begin{array}{lcl}
a&=&(c,b),\\
b&=&(b,c),\\
c&=&(d,d)\sigma,\\
d&=&(a,a)\sigma.
\end{array}
\end{equation*}

\begin{theorem}[]
The group generated by the automaton $\B_4$ is the free product of
$4$ copies of cyclic group of order 2.
\end{theorem}

For any $n>4$ let $\B^{(n)}$ be the $n$-state automaton defined by
the wreath recursion
\begin{equation*}
\begin{array}{lcl}
a_n&=&(c_n,b_n),\\
b_n&=&(b_n,c_n),\\
c_n&=&(q_{n1},q_{n1})\sigma,\\
q_{n,i}&=&(q_{n,i+1},q_{n,i+1})\sigma_{n,i},i=1,\ldots,n-5,\\
q_{n,n-4}&=&(d_n,d_n)\sigma_{n,n-4},\\
d_n&=&(a_n,a_n)\sigma,\\
\end{array}
\end{equation*}
where $\sigma_{n,i}\in\Sym(\{1,2\})$ are chosen arbitrarily.

\begin{theorem}
The group generated by the automaton $\B^{(n)}$ is the free product
of $n$ copies of cyclic group of order 2.
\end{theorem}

The Bellaterra automaton $\B_3$ with 3 states is closely related to
a 3-state automaton introduced by Aleshin~\cite{aleshin:free} in the
beginning of 1980's. Namely, the only difference is in the output
function, whose value is always opposite to the one of the Aleshin
automaton. It was proved by M.~Vorobets and
Y.~Vorobets~\cite{vorobets:aleshin} that the group generated by this
automaton is a free group of rank 3.

The series of automata from~\cite{vorobets:series_free} and the
family of automata constructed in~\cite{steinberg_vv:series_free}
generating the free products of groups of order 2 have counterpart
series and family of automata generalizing the Aleshin 3-state
automaton and generating the free groups. In this paper the proofs
for free products are simpler, but, as a downside, our technique
does not give an answer to the question whether changing the output
function of automata in our family to the opposite value produces
the automata generating the free groups.

We can not miss mentioning the amusing fact that among all 190 types
of 3-state automata acting on the 2-letter alphabet, that are not
symmetric to each other the only automaton generating a free
non-abelian group is the first Aleshin automaton
(see~\cite{bondarenko_gkmnss:full_clas32}). Thus pinpointing by
Aleshin this unique automaton a quarter of a century ago should
really be highly valued.

The structure of the paper is as follows. All necessary definitions
are given in Section~\ref{sec_prelim}. The automaton generating the
free product of $4$ cyclic groups of order 2 is studied in
Section~\ref{sec_base_case}. In Section~\ref{sec_series} the family
of automata generating the free products of groups of order 2 is
considered.

\section{Preliminaries}
\label{sec_prelim}

Let $X$ be a finite set of cardinality $d$. By $X^*$ we denote the
free monoid generated by $X$, which consists of finite words over
$X$. This monoid can be naturally endowed with a structure of a
rooted $d$-ary tree by declaring that $v$ is adjacent to $vx$ for
any $v\in X^*$ and $x\in X$. The empty word corresponds to the root
of the tree and $X^n$ corresponds to the $n$-th level of the tree.
We will be interested in the groups of automorphisms and semigroups
of homomorphisms of $X^*$. Any such homomorphism can be defined via
the notion of initial automaton.

\begin{defin}
A \emph{Mealy automaton} (or simply \emph{automaton}) is a tuple
$(Q,X,\pi,\lambda)$, where $Q$ is a set (a set of states), $X$ is a
finite alphabet, $\pi\colon Q\times X\to Q$ is a transition function
and $\lambda\colon Q\times X\to X$ is an output function. If the set
of states $Q$ is finite the automaton is called \emph{finite}. If
for every state $q\in Q$ the output function $\lambda(q,x)$ induces
a permutation of $X$, the automaton $\A$ is called invertible.
Selecting a state $q\in Q$ produces an \emph{initial automaton}
$\A_q$.
\end{defin}

\begin{figure}
\begin{center}
\epsfig{file=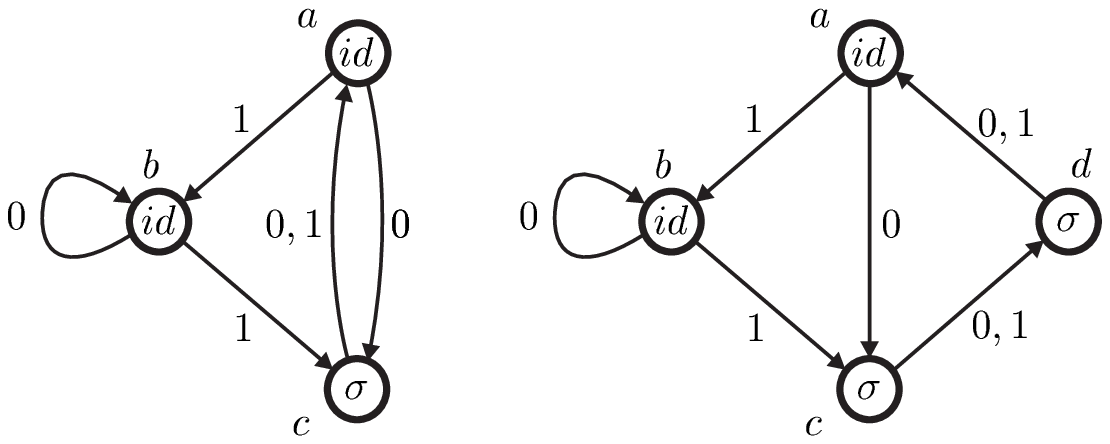} \caption{Bellaterra automata
$\B_3$ and $\B_4$\label{fig_bel}}
\end{center}
\end{figure}

Automata are often represented by the \emph{Moore diagrams}. The
Moore diagram of an automaton $\A=(Q,X,\pi,\lambda)$ is a directed
graph in which the vertices are the states from $Q$ and the edges
have form $q\stackrel{x|\lambda(q,x)}{\longrightarrow}\pi(q,x)$ for
$q\in Q$ and $x\in X$. If the automaton is invertible, then it is
common to label vertices of the Moore diagram by the permutation
$\lambda(q,\cdot)$ and leave just first components from the labels
of the edges. An example of Moore diagram is shown in
Figure~\ref{fig_bel}.

Any initial automaton induces a homomorphism of $X^*$. Given a word
$v=x_1x_2x_3\ldots x_n\in X^*$ it scans its first letter $x_1$ and
outputs $\lambda(x_1)$. The rest of the word is handled in a similar
fashion by the initial automaton $\A_{\pi(x_1)}$. Formally speaking,
the functions $\pi$ and $\lambda$ can be extended to $\pi\colon
Q\times X^*\to Q$ and $\lambda\colon  Q\times X^*\to X^*$ via
\[\begin{array}{l}
\pi(q,x_1x_2\ldots x_n)=\pi(\pi(q,x_1),x_2x_3\ldots x_n),\\
\lambda(q,x_1x_2\ldots x_n)=\lambda(q,x_1)\lambda(\pi(q,x_1),x_2x_3\ldots x_n).\\
\end{array}
\]

By construction any initial automaton acts on $X^*$ as a
homomorphism. In case of invertible automaton it acts as an
automorphism.

\begin{defin}
The semigroup (group) generated by all states of automaton $\A$ is
called an \emph{automaton semigroup} (\emph{automaton group}) and
denoted by $\mathds S(\A)$ (respectively $\mathds G(\A)$).
\end{defin}

Another popular name for automaton groups and semigroups is
self-similar groups and semigroups
(see~\cite{nekrash:self-similar}).

Conversely, any homomorphism of $X^*$ can be encoded by the action
of an initial automaton. In order to show this we need a notion of a
\emph{section} of a homomorphism at a vertex of the tree. Let $g$ be
a homomorphism of the tree $X^*$ and $x\in X$. Then for any $v\in
X^*$ we have
\[g(xv)=g(x)v'\]
for some $v'\in X^*$. Then the map $g|_x\colon X^*\to X^*$ given by
\[g|_x(v)=v'\]
defines a homomorphism of $X^*$ and is called the \emph{section} of
$g$ at vertex $x$. Furthermore,  for any $x_1x_2\ldots x_n\in X^*$
we define \[g|_{x_1x_2\ldots x_n}=g|_{x_1}|_{x_2}\ldots|_{x_n}.\]

Given a homomorphism $g$ of $X^*$ we construct an initial automaton
$\A(g)$ whose action on $X^*$ coincides with that of $g$ as follows.
The set of states of $\A(g)$ is the set $\{g|_v\colon  v\in X^*\}$
of different sections of $g$ at the vertices of the tree. The
transition and output functions are defined by
\[\begin{array}{l}
\pi(g|_v,x)=g|_{vx},\\
\lambda(g|_v,x)=g|_v(x).
\end{array}\]

Throughout the paper we will use the following convention. If $g$
and $h$ are the elements of some (semi)group acting on set $A$ and
$a\in A$, then
\begin{equation}
\label{eqn_conv}
gh(a)=h(g(a)).
\end{equation}

Taking into account convention~\eqref{eqn_conv} one can compute
sections of any element of an automaton semigroup as follows. If
$g=g_1g_2\cdots g_n$ and $v\in X^*$, then

\begin{equation}
\label{eqn_sections} g|_v=g_1|_v\cdot g_2|_{g_1(v)}\cdots
g_n|_{g_1g_2\cdots g_{n-1}(v)}.
\end{equation}

For any automaton group $G$ there is a natural embedding
\[G\hookrightarrow G \wr \Sym(X)\]
defined by
\[G\ni g\mapsto (g_1,g_2,\ldots,g_d)\lambda(g)\in G\wr \Sym(X),\]
where $g_1,g_2,\ldots,g_d$ are the sections of $g$ at the vertices
of the first level, and $\lambda(g)$ is a permutation of $X$ induced
by the action of $g$ on the first level of the tree.

The above embedding is convenient in computations involving the
sections of automorphisms, as well as for defining automaton groups.
For example, the group $\mathds G(\B_3)$ generated by automaton in
Figure~\ref{fig_bel}, where $\sigma=(1,2)$ denotes the nontrivial
element of $\Sym(\{1,2\})$, can be defined as
\begin{equation}
\label{eqn_bella_def}
\begin{array}{lcl}
a&=&(c,b),\\
b&=&(b,c),\\
c&=&(a,a)\sigma.\\
\end{array}
\end{equation}
The latter definition is sometimes called the \emph{wreath
recursion} defining the group.

For any finite automaton one can construct a \emph{dual automaton}
defined by switching the states and the alphabet as well as
switching the transition and the output functions.

\begin{defin}
Given a finite automaton $\A=(Q,X,\pi,\lambda)$ its \emph{dual
automaton} $\hat\A$ is a finite automaton
$(X,Q,\hat\lambda,\hat\pi)$, where
\[\begin{array}{l}
\hat\lambda(x,q)=\lambda(q,x),\\
\hat\pi(x,q)=\pi(q,x)
\end{array}\]
for any $x\in X$ and $q\in Q$.
\end{defin}

Note that the dual of the dual of an automaton $\A$ coincides with
$\A$. The semigroup $\mathds S(\hat\A)$ generated by dual automaton
$\hat\A$ of automaton $\A$ acts on the free monoid $Q^*$. This
action induces the action on $\mathds S(\A)$. Similarly, $\mathds
S(\A)$ acts on $\mathds S(\hat\A)$.

\begin{defin}
For an automaton semigroup $G$ generated by automaton $\A$ the
\emph{dual semigroup} $\hat G$ to $G$ is a semigroup generated by a
dual automaton $\hat\A$.
\end{defin}

A particularly important class of automata is the class of
bireversible automata.

\begin{defin}
An automaton $\A$ is called \emph{bireversible} if it is invertible,
its dual is invertible, and the dual to $\A^{-1}$ are invertible.
\end{defin}

In particular, for any group generated by a bireversible automaton
$\A$ one can consider a dual group generated by the dual auomtaton
$\hat\A$.

The following proposition is proved in~\cite{vorobets:aleshin} and
is proved by induction on level. With a slight abuse of notations we
will denote by the same symbol the element of a free monoid and its
image under canonical epimorphism onto corresponding semigroup.

\begin{prop}
\label{prop_dual_sections}
Let $G$ be an automaton semigroup acting on $X^*$ and generated by
the finite set $S$.  And let $\hat G$ be a dual semigroup to $G$
acting on $S^*$. Then for any $g\in G$ and $v\in X^*$ we have
$g|_v=v(g)$ in $G$. Similarly, for any $g\in S^*$ and $v\in \hat G$,
$v|_g=g(v)$ in $\hat G$.
\end{prop}

\section{Automaton generating $C_2\ast C_2\ast C_2\ast C_2$}
\label{sec_base_case}

Consider the group $\G$ generated by the following 4-state automaton
$\B_4$, whose transition and output functions are given by wreath
recursion (its Moore diagram is shown in the right half of
Figure~\ref{fig_bel})
\begin{equation}
\label{eqn_autom_def}
\begin{array}{lcl}
a&=&(c,b),\\
b&=&(b,c),\\
c&=&(d,d)\sigma,\\
d&=&(a,a)\sigma.
\end{array}
\end{equation}

The next Theorem is the main result of this section.

\begin{theorem}
\label{thm_main}
Group $\G$ is isomorphic to $C_2\ast C_2\ast C_2\ast C_2$.
\end{theorem}

The proof of this theorem is split into a number of lemmas below.

First, we note that the automaton $\B_4$ is bireversible. The dual
group $\Gamma$ to $\G$ is generated by the following automaton
\begin{equation}
\label{eqn_autom_dual_def}
\begin{array}{lcl}
\zero&=&(\zero,\zero,\one,\one)(a\,c\,d),\\
\one&=&(\one,\one,\zero,\zero)(a\,b\,c\,d).\\
\end{array}
\end{equation}

\begin{figure}[h]
\begin{center}
\epsfig{file=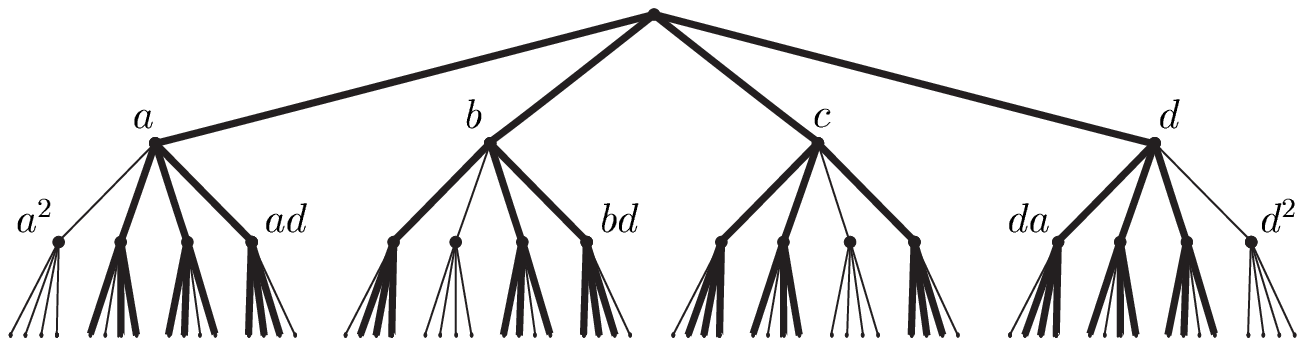, height=75pt} \caption{Trees
$T$ and $\hat T$\label{fig_trees}}
\end{center}
\end{figure}

This group acts on a rooted 4-ary tree $T$ whose vertices are
labelled by the words over $\{a,b,c,d\}$. Since $a^2=(c^2,b^2)$,
$b^2=(b^2,c^2)$, $c^2=(d^2,d^2)$ and $d^2=(a^2,a^2)$ we get that
$a^2=b^2=c^2=d^2=1$ in $\Gamma$ and the image of any word containing
any of $a^2$, $b^2$, $c^2$ or $d^2$ under any element of $\Gamma$
will also contain one of these subwords. Therefore there is an
invariant under $\Gamma$ subtree $\hat T$ of $T$ consisting of all
words over $\{a,b,c,d\}$ that do not have $a^2$, $b^2$, $c^2$ and
$d^2$ as subwords. The root of $\hat T$ has 4 descendants and all
the other vertices in $\hat T$ have three (see
Figure~\ref{fig_trees}, where subtree $\hat T$ is drawn with bold
edges).

The following simple proposition was obtained independently by Z.~\v
Suni\'c (private communication) and the proof is implicitly
contained in the book of Nekrashevych~\cite{nekrash:self-similar}.

\begin{prop}
\label{prop_dual_finite}
Let $G$ be any semigroup generated by a finite automaton and $\hat
G$ be its dual semigroup. Then $G$ is finite if and only if $\hat G$
is finite.
\end{prop}
\begin{proof}
Since dual of the dual of the automaton generating $G$ coincides
with this automaton, it is enough to show the implication in one
direction.

Suppose $G$ is finite. For any element $v\in \hat G$ and any vertex
$g$ of tree the semigroup $\hat G$ acts on, we have $v|_g=g(v)$ in
$\hat G$ by Proposition~\ref{prop_dual_sections}. Therefore the
number of different sections of $v$ is bounded by the size of $G$.
But there are only finitely many different automata with a fixed
number of states. Thus $\hat G$ is finite.
\end{proof}

\begin{lemma}
\label{lem_G_infinite}
The group $\G$ is infinite.
\end{lemma}
\begin{proof}
The lemma follows from the fact that the group acts transitively on
each level of the tree. To prove this we first observe that the
group $G/\Stg(2)$ is cyclic of order 4 and the portrait of depth 2
of every element of $G$ (rooted binary tree of depth $2$, where each
vertex is labelled by the permutation induced by this element at
this vertex) must coincide with one of the listed in
Figure~\ref{fig_portraits}.

\begin{figure}[h]
\begin{center}
\epsfig{file=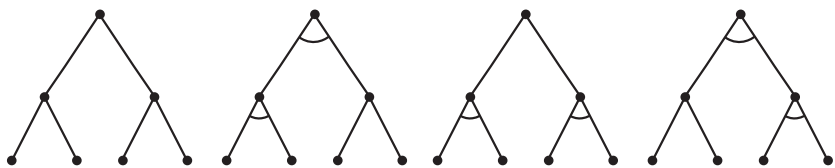, height=60pt} \caption{Possible portraits
of elements of $G$ of depth 2\label{fig_portraits}}
\end{center}
\end{figure}

It is proved in~\cite{gawron_ns:conjugation} that an automorphism
$g$ of the rooted binary tree acts level transitively if and only if
on each level the number of sections of $g$ at the vertices of this
level acting nontrivially on the first level, is odd.

By induction on level it follows that each element $g$ of $G$ acting
nontrivially on the first level acts spherically transitively.
Indeed, if the number of sections of $g$ at the vertices of the
$k$-th level acting nontrivially on the first level (the number of
``switches'' on the $k$-th level) is odd, then each of these
sections will produce exactly one switch on the $(k+1)$-st level,
while the sections acting trivially on the first level will produce
either none or two switches on the $(k+1)$-st level. Thus, the total
number of switches on the $(k+1)$-st level will be odd as well.
\end{proof}

The direct corollary of Proposition~\ref{prop_dual_finite} and
Lemma~\ref{lem_G_infinite} is
\begin{cor}
\label{cor_Gamma_infinite}
The group $\Gamma$ is infinite.
\end{cor}

\begin{cor}
\label{cor_stabs_diff}
The stabilizers of levels of $T$ in $\Gamma$ are pairwise different.
\end{cor}
\begin{proof}
Since $\Gamma$ is infinite by Corollary~\ref{cor_Gamma_infinite} and
all stabilizers of levels are finite index subgroups in $\Gamma$,
they are all infinite. Let $g\in\Stga(n)$ be arbitrary and
nontrivial and let $m\geq n+1$ be the smallest level on which $g$
acts nontrivially. Then there exists a vertex $v=x_1x_2\ldots
x_{m-1}$ of the tree, such that $g|_v$ acts nontrivially on the
first level. Then $g|_{x_1x_2\ldots x_{m-n-1}}\in
\Stga(n)\setminus\Stga(n+1)$.
\end{proof}

\begin{lemma}
\label{lem_stabs_equal}
Let $\hat T_n$ be the subtree of $\hat T$ consisting of the first
$n$ levels. Then $\Stga(n)=\Stga(\hat T_n)$.
\end{lemma}

\begin{proof}
Since the leaves of $\hat T_n$ are vertices of the $n$-th level of
$T$ we have $\Stga(n)\subset\Stga(T_n)$.

Suppose $v\in\Stga(\hat T_n)\setminus\Stga(n)$. Then there is a
vertex $g$ of the $n$-th level which is not in $\hat T$ and is not
fixed under $v$. Since $v$ fixes $T_n$ it follows that $g=ftth$ and
$v(g)=ftth'$ for some $f,h,h'\in G$ and $t\in\{a,b,c,d\}$. Then
\[v(fh)=v(f)v|_f(h)=f\cdot(v|_f)|_{tt}(h)=fv|_{ftt}(h)=fh'.\]
The second equality above holds since for any $t\in\{a,b,c,d\}$ we
have $t^2=1$ and thus for any $w\in\Gamma$ by
Proposition~\ref{prop_dual_sections}, $w|_{tt}=(tt)(w)=w$ in
$\Gamma$ and for any word $h\in T$ we have $w|_{tt}(h)=w(h)$.

We can repeat this procedure until we get an element of $\hat T_n$
not fixed under the action of $v$, obtaining contradiction. Thus
$v\in\Stga(n)\setminus\Stga(T_n)$.
\end{proof}

The next statement follows directly from
Corollary~\ref{cor_stabs_diff} and Lemma~\ref{lem_stabs_equal}.

\begin{cor}
\label{cor_stabs_Tn_prime}
For any $n\geq1$ there is an element in $\Gamma$ fixing $\hat T_n$
but moving some vertex in $\hat T_{n+1}$.
\end{cor}

\begin{lemma}
\label{lem_even_perms}
The sections of any element of $\Stga(n)$ at the vertices of the
$n$-th level act on the first level by even permutations.
\end{lemma}

\begin{proof}
By self-similarity it is enough to prove the Lemma for $n=1$. The
claim follows from the fact that $|\Stga(1)/\Stga(2)|=3^3$ (this was
computed using~\cite{muntyan_s:automgrp}). Therefore the sections of
any element of $\Stga(1)$ at the vertices of the first level act on
the first level by permutations, whose order is a power of 3, which
are either cycles of length 3 or the trivial permutation. All these
are even permutations.

Below we provide a proof that does not rely on computer
computations. This proof is also important because it introduces
certain notation that will be used later in
Section~\ref{sec_series}.

First we show that if $v\in\Gamma$ fixes vertex $d$, then the
parities of the actions of $v$ and $v|_d$ on the first level
coincide. For this purpose we introduce a new generating set in
$\Gamma$. For any $x\in\Sym(\{a,b,c,d\})$ denote by $\bar x$ the
automorphism of $T$ defined by
\[\bar x=(\bar x,\bar x,\bar x,\bar x)x.\]
The portrait of $\bar x$ has $x$ at each vertex of the tree. Since
$(\zero\one^{-1})^2=(\zero^{-1}\one)^2=1$ we obtain

\begin{equation}
\label{eqn_ab}
\zero\one^{-1}=(\zero\one^{-1},\zero\one^{-1},\zero\one^{-1},\zero\one^{-1})(a\,b)=\overline{(a\,b)}
\end{equation}
and
\begin{equation}
\label{eqn_bc}
\zero^{-1}\one=(\zero^{-1}\one,\zero^{-1}\one,\zero^{-1}\one,\zero^{-1}\one)(b\,c)=\overline{(b\,c)}.
\end{equation}
This shows that $\overline{(a\,c)}\in\Gamma$. If we denote
\begin{equation}
\label{eqn_alpha_beta}
\begin{array}{lclr}
\alpha=\one\cdot\overline{(a\,c)}&=&(\alpha,\alpha,\beta,\beta)&(a\,b)(c\,d),\\
\beta
=\zero\cdot\overline{(a\,c)}&=&(\beta,\beta,\alpha,\alpha)&(c\,d),
\end{array}
\end{equation}
then $\alpha^2=\beta^2=1$ and, taking into account that
$\beta\alpha^{-1}=\overline{(a\,b)}$,
\[\Gamma=\langle\alpha,\beta,\overline{(b\,c)}\rangle.\]

Suppose now $v\in\Gamma$ is an arbitrary element fixing vertex $d$.
Represent $v$ as a word over $\{\alpha,\beta,\overline{(b\,c)}\}$
\[v=v_1v_2\cdots v_k,\]
then by~\eqref{eqn_sections}
\[v|_d=v_1|_d\cdot v_2|_{v_1(d)}\cdots v_k|_{v_1v_2\cdots v_{k-1}(d)}.\]
The parity of the action of $v_i$ on the first level is different
from the one of $v_i|_{v_{1}v_{2}\cdots v_{i-1}(d)}$ only in case
$v_i$ is $\alpha$ or $\beta$ and $v_{1}v_{2}\cdots v_{i-1}(d)=c$ or
$v_{1}v_{2}\cdots v_{i-1}(d)=d$.

Note that in this situation if $v_{1}v_2\cdots v_{i-1}(d)=c$ then
$v_1v_{2}\cdots v_i(d)=d$, and if $v_{1}v_2\cdots v_{i-1}(d)=d$ then
$v_{1}v_2\cdots v_{i}(d)=c$. The converse is also true in the
following sense: if $v_{1}v_{2}\cdots v_{i-1}(d)\neq d$ and
$v_1v_{2}\cdots v_i(d)=d$ then $v_{1}v_{2}\cdots v_{i-1}(d)=c$ and
$v_i$ is either $\alpha$ or $\beta$, and if $v_{1}v_{2}\cdots
v_{i-1}(d)=d$ and $v_1v_{2}\cdots v_i(d)\neq d$, then
$v_1v_{2}\cdots v_i(d)=c$ and $v_i$ is either $\alpha$ or $\beta$.
In other words, the parity of the action of $v_i$ on the first level
is different from the one of $v_i|_{v_{1}v_2\cdots v_{i-1}(d)}$
exactly when there is a change from $d$ to anything else or from
something to $d$ in the sequence $\{d, v_1(d),\ldots,v_1v_2\cdots
v_k(d)\}$. But since $v_1v_2\cdots v_k(d)=v(d)=d$, there must be an
even number of such changes. Hence, the parity is different in even
number of places and the parities of the actions of $v$ and $v|_d$
on the first level coincide.

By the above for any $g=(g|_a,g|_b,g|_c,g|_d)\in\Stga(1)$ the parity
of the action of $g|_d$ on the first level is even. Furthermore, the
conjugate $g^\beta=\beta^{-1}g\beta$ has decomposition
% a^b=b^-1ab
\[g^{\beta}=(\ast,\ast,\ast,(g|_c)^\alpha)\in \Stga(1),\]
which implies that $(g|_c)^\alpha$ and, hence, $g|_c$ acts on the
first level by an even permutation. Finally,
\[\begin{array}{l}
g^{\overline{(a\,c)}}=(\ast,\ast,(g|_a)^{\overline{(a\,c)}},\ast)\in
\Stga(1),\\
g^{\overline{(b\,c)}}=(\ast,\ast,(g|_b)^{\overline{(b\,c)}},\ast)\in
\Stga(1).
\end{array}\]
This shows that all sections of $g$ at the vertices of the first
level act on the first level by even permutations.
\end{proof}

\begin{lemma}
\label{lem_transit}
The group $\Gamma$ acts transitively on the levels of $\hat T$.
\end{lemma}

\begin{proof}
We proceed by induction on levels. The transitivity on the first
level is clear. Assume $\Gamma$ acts transitively on the $n$-th
level of $\hat T$. By Corollary~\ref{cor_stabs_Tn_prime} there is an
element $v\in\Gamma$ that fixes $\hat T_n$ and acts nontrivially on
$\hat T_{n+1}$. This means that there is a vertex $g\in \hat T_n$
such that $v(g)=g$ and $v|_g$ acts nontrivially on the first level.
By Lemma~\ref{lem_even_perms} the permutation induced by $v|_g$ on
the first level is even, which implies that it is a cycle of length
3. Thus $v|_g$ acts transitively on the first level of the tree.

Without loss of generality assume that $g$ ends with $d$. By
induction assumption, for any vertex $h_1h_2\ldots h_{n+1}$ of $\hat
T_{n+1}$ there is an element $w\in\Gamma$ that moves $g$ to
$h_1h_2\ldots h_n$. Then $v^kw$, where $k$ is 0, 1 or 2 will move
$ga$ to $h_1h_2\ldots h_{n+1}$. Thus, $\Gamma$ acts transitively on
$\hat T_{n+1}$.
\end{proof}

Finally, we have all ingredients for the proof of
Theorem~\ref{thm_main}.

\begin{proof}[Proof of Theorem~\ref{thm_main}]
For every $n\geq1$ there is a nontrivial element $h\in\G$ that
belongs to the $n$-th level of $\hat T$
($h=(ab)^{\frac{n-1}2}c\neq1$ for an odd $n$ and $h=(ab)^{\frac
n2-1}ac\neq 1$ for an even $n$). By Lemma~\ref{lem_transit} the
group $\Gamma$ acts transitively on each level of $\hat T$.
Therefore for any word $g$ from the $n$-th level of $\hat T$ (which
is a word of length $n$ without double letters) there exists
$v\in\Gamma$ such that
\[g|_v=v(g)=h\neq1.\]
Thus there are no relations in $\G$ except $a^2=b^2=c^2=d^2=1$.
\end{proof}

\section{Family of automata generating the free products of $C_2$}
\label{sec_series}

\begin{figure}
\begin{center}
\epsfig{file=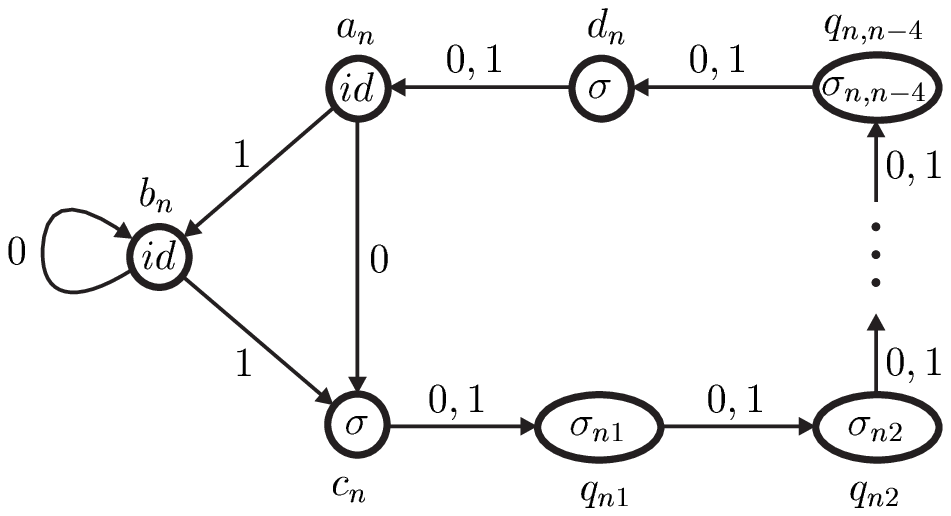} \caption{Bellaterra
automaton $\B^{(n)}$\label{fig_bel_n}}
\end{center}
\end{figure}

Let us define a family of automata obtained from the automaton
$\B_4$ by inserting new states on the path from $c$ to $d$. Namely,
for every integer $n>4$ and any permutations
$\sigma_{n,i}\in\Sym(\{1,2\})$, $i=1,\ldots,n-4$ consider an
automaton with $n$ states
$a_n,b_n,c_n,q_{n1},q_{n2},\ldots,q_{n,n-4},d_n$ whose transition
and output functions are given via the wreath recursion
\begin{equation}
\label{eqn_series_def}
\begin{array}{lcl}
a_n&=&(c_n,b_n),\\
b_n&=&(b_n,c_n),\\
c_n&=&(q_{n1},q_{n1})\sigma,\\
q_{n,i}&=&(q_{n,i+1},q_{n,i+1})\sigma_{n,i},i=1,\ldots,n-5,\\
q_{n,n-4}&=&(d_n,d_n)\sigma_{n,n-4},\\
d_n&=&(a_n,a_n)\sigma.\\
\end{array}
\end{equation}
With a slight abuse of notation we denote this automaton by
$\B^{(n)}$ regardless of the choice of permutations $\sigma_{n,i}$.
The Moore diagram of $\B^{(n)}$ is shown in Figure~\ref{fig_bel_n}.

From the wreath recursion it is easy to observe that the generators
of $\B^{(n)}$ are involutions.

This section is devoted to the proof of the following theorem.

\begin{theorem}
\label{thm_general_case}
The group $\G^{(n)}$ generated by automaton $\B^{(n)}$ is isomorphic
to the free product of $n$ copies of the cyclic group of order $2$.
\end{theorem}

The proof relies on the results of Section~\ref{sec_base_case}.  The
approach is similar.  We prove that the dual automaton acts
transitively on the invariant subtree consisting of words without
double letters.  This yields the structure of the free product in
the group $\G^{(n)}$.

Note that the automaton $\B^{(n)}$ is bireversible so that the dual
group $\Gamma^{(n)}$ to $\G^{(n)}$ is well defined. The group
$\Gamma^{(n)}$ is generated by automaton acting on the rooted
$n$-ary tree $\Tn$ as follows

\begin{equation}
\label{eqn_autom_dual_def_series}
\begin{array}{l}
\zero_n=(\zero_n,\zero_n,\one_n,\mathds K_{n1},\ldots,\mathds K_{n,n-4},\one_n)(a_n\,c_n\,q_{n1}\ldots q_{n,n-4}\,d_n),\\
\one_n=(\one_n,\one_n,\zero_n,\mathds L_{n1},\ldots,\mathds L_{n,n-4},\zero_n)(a_n\,b_n\,c_n\,q_{n1}\ldots q_{n,n-4}\,d_n),\\
\end{array}
\end{equation}
where $\mathds K_{n,i}=\zero_n$ and $\mathds L_{n,i}=\one_n$ if
$\sigma_{n,i}$ is a trivial permutation, and $\mathds
K_{n,i}=\one_n$ and $\mathds L_{n,i}=\zero_n$ otherwise.

Consider a subtree $\hat T^{(n)}$ of $\Tn$ consisting of all words
over the alphabet
$Y^{(n)}=\{a_n,b_n,c_n,q_{n1},q_{n2},\ldots,q_{n,n-4},d_n\}$ without
double letters. The root of $\hat T^{(n)}$ has $n$ descendants and
all other vertices have $n-1$. This subtree is invariant under the
action of $\Gamman$.

Similarly to~\eqref{eqn_ab} and~\eqref{eqn_bc} we get that
$\zero_n\one_n^{-1}=\overline{(a_n\,b_n)}$ and
$\zero_n^{-1}\one_n=\overline{(b_n\,c_n)}$. Similarly
to~\eqref{eqn_alpha_beta} we define transformations
$\alpha_n=\one_n\cdot\overline{(a_n\,c_n)}$ and $\beta_n
=\zero_n\cdot\overline{(a_n\,c_n)}$ for which we have
\begin{equation}
\label{eqn_alpha_beta_n}
\begin{array}{lr}
\alpha_n=(\alpha_n,\alpha_n,\beta_n,\gamma_{n1},\ldots,\gamma_{n,n-4},\beta_n)&(a_n\,b_n)(c_n\,q_{n1}\ldots q_{n,n-4}\,d_n),\\
\beta_n=(\beta_n,\beta_n,\alpha_n,\delta_{n1},\ldots,\delta_{n,n-4},\alpha_n)&(c_n\,q_{n1}\ldots
q_{n,n-4}\,d_n),
\end{array}
\end{equation}
where $\gamma_{n,i}=\alpha_n$ and $\delta_{n,i}=\beta_n$ if
$\sigma_{n,i}$ is a trivial permutation, and $\gamma_{n,i}=\beta_n$
and $\delta_{n,i}=\alpha_n$ otherwise.

Since $\alpha_n^{-1}\beta_n=\overline{(a_n\,b_n)}$ we get a new
generating set for $\Gamma_n$,
\[\Gamman=\langle\alpha_n,\beta_n,\overline{(b_n\,c_n)}\rangle.\]

The following lemma establishes a relation between the actions of
the groups $\Gamma$ and $\Gamman$. We consider the tree $\hat T$
naturally embedded in the tree $\hat T^{(n)}$ via a homomorphism of
monoids induced by $a\mapsto a_n$, $b\mapsto b_n$, $c\mapsto c_n$,
$c\mapsto c_n$, $d\mapsto d_n$. Then the group $\Gamma$ acts also on
$\hat T^{(n)}$ (the action on the letters not in the image of $\hat
T$ is defined to be trivial).

\begin{lemma}
\label{lem_actions}
For any $v\in\Gamma$ there exists $v'\in\Gamman$ with the following
property. For any word $g$ over $\{a_n,b_n,c_n\}$ such that $v(g)$
is also a word over $\{a_n,b_n,c_n\}$, we have $v(g)=v'(g)$.
\end{lemma}

\begin{proof}
Let $x_1x_2\ldots x_k$ be the word over
$\{\alpha,\beta,\overline{(b_n\,c_n)}\}$ representing $v$. Define
$y_i\in\{\alpha_n,\beta_n,\overline{(b_n\,c_n)}\}$ by the following
rule. If $x_i=\overline{(b_n\,c_n)}$, then put $y_i=x_i$. In the
case $x_i=\alpha$ (resp. $x_i=\beta$) compute the total number of
$\alpha$ and $\beta$  among $x_1,x_2,\ldots,x_{i-1}$. If this number
is even, then define $y_i=\alpha_n$ (resp. $y_i=\beta_n$).
Otherwise, put $y_i=\alpha_n^{-1}$ (resp. $y_i=\beta_n^{-1}$).

Now let $g$ be any word over $\{a_n,b_n,c_n\}$. We will show by
induction on $i$ that $y_1y_2\ldots y_i(g)$ is obtained from
$x_1x_2\ldots x_i(g)$ by replacing all occurrences of $d_n$ by
$q_{n1}$ when the total number of $\alpha$ and $\beta$ among
$x_1,x_2,\ldots,x_{i-1}$ is odd, and coincides with $x_1x_2\ldots
x_i(g)$ otherwise.

The claim holds trivially for $i=0$. Let us prove the induction
step. First of all, if $x_{i+1}=y_{i+1}=\overline{(b_n\,c_n)}$ then
the relation between $y_1y_2\ldots y_{i+1}(g)$ and $x_1x_2\ldots
x_{i+1}(g)$ is the same as between $y_1y_2\ldots y_{i}(g)$ and
$x_1x_2\ldots x_{i}(g)$. This is because $\overline{(b_n\,c_n)}$
fixes letters $d_n$ and $q_{n,1}$. Hence we can assume that
$x_{i+1}=\alpha$ or $x_{i+1}=\beta$.

Suppose first that there is an odd number of $\alpha$ and $\beta$
among $x_1,x_2,\ldots,x_i$. By induction assumption $y_1y_2\ldots
y_i(g)$ is obtained from $x_1x_2\ldots x_i(g)$ by replacing all
occurrences of $d_n$ by $q_{n1}$ and, in particular, is a word over
$\{a_n,b_n,c_n,q_{n1}\}$. If $x_{i+1}=\alpha$ ($x_{i+1}=\beta$),
then by construction $y_{i+1}=\alpha_n^{-1}$ (respectively
$y_{i+1}=\beta_n^{-1}$), for which we have

\begin{equation}
\label{eqn_alpha_beta_n_inv}
\begin{array}{l}
\alpha_n^{-1}=(\alpha_n^{-1},\alpha_n^{-1},\beta_n^{-1},\beta_n^{-1},\gamma_{n1}^{-1},\ldots,\gamma_{n,n-4}^{-1})(a_n\,b_n)(c_n\,d_n\ldots q_{n1}),\\
\beta_n^{-1}=(\alpha_n^{-1},\alpha_n^{-1},\beta_n^{-1},\beta_n^{-1},\delta_{n1}^{-1},\ldots,\delta_{n,n-4}^{-1})\,\phantom{(a_n\,b_n)}(c_n\,d_n\ldots
q_{n1}).
\end{array}
\end{equation}

Therefore the images of $y_1y_2\ldots y_{i}(g)$ under the actions of
$\alpha_n^{-1}$ and $\beta_n^{-1}$ coincide with the images of
$x_1x_2\ldots x_i(g)$ under the actions of $\alpha$ and $\beta$
correspondingly. Thus, $y_1y_2\ldots y_{i+1}(g)=x_1x_2\ldots
x_{i+1}(g)$, which is exactly what we need since the number of
$\alpha$ and $\beta$ among $x_1,x_2,\ldots,x_{i+1}$ is even.

In case of even number of occurrences of $\alpha$ and $\beta$ among
$x_1,x_2,\ldots,x_{i}$ by induction assumption $y_1y_2\ldots y_i(g)$
coincides with $x_1x_2\ldots x_i(g)$. In particular, it is a word
over $\{a_n,b_n,c_n,d_n\}$. Also by construction $y_{n+1}=\alpha_n$
or $y_{n+1}=\beta_n$.

It follows from~\eqref{eqn_alpha_beta_n} that $y_{i+1}$ acts on the
letters of $y_1y_2\ldots y_i(g)$ exactly as $x_{i+1}$, except it
everywhere moves $c_n$ to $q_{n1}$, instead of moving it to $d_n$.
Therefore, the resulting word $y_1y_2\ldots y_{i+1}(g)$ can be
obtained from $x_1x_2\ldots x_{i+1}(g)$ by changing all occurrences
of $d_n$ by $q_{n1}$. This agrees with the fact that the total
number of $\alpha$ and $\beta$ among $x_1,x_2,\ldots,x_{i+1}$ is
odd.

Finally, to finish the proof of the lemma, it is enough to put
$v'=y_1y_2\ldots y_k$ and note that if $v(g)$ is a word over
$\{a_n,b_n,c_n\}$, then $v'(g)$ must coincide with $v(g)$ regardless
of the number of $\alpha$ and $\beta$ in the word representing $v$.
\end{proof}

\begin{lemma}
\label{lem_gamman_trans}
The group  $\Gamman$ acts transitively on the levels of $\hat
T^{(n)}$.
\end{lemma}
\begin{proof}
We proceed by induction on levels. Obviously, $\Gamman$ acts
transitively on the first level. Suppose it acts transitively on
level $m$. We will show that any vertex of the $(m+1)$-st level can
be moved to the vertex $a_nb_na_nb_n\ldots b_na_n$ or
$a_nb_na_nb_n\ldots a_nb_n$ (depending on the parity of $m$).

Let $g$ be the vertex of the $(m+1)$-st level of $\hat T^{(n)}$.
Then $g=ht$, where $h$ is the vertex of the $m$-th level and $t\in
Y^{(n)}$. For definiteness let us assume that $m$ is even. By
induction assumption there is $v\in\Gamman$ that moves $h$ to
$a_nb_na_nb_n\ldots b_n$. Then
\[v(g)=a_nb_na_nb_n\ldots a_nb_nt'\]
for some $t'\in Y^{(n)}$. Since $\beta_n$ fixes $a_nb_na_nb_n\ldots
a_nb_n$ and $\beta_n|_{a_nb_na_nb_n\ldots a_nb_n}=\beta_n$, after
applying, if necessary, a power of $\beta_n$ we can assume that
$t'\in\{a_n,b_n,c_n\}$. Now we invoke the transitivity of the group
$\Gamma$ on $\hat T$. By Lemma~\ref{lem_transit} there is
$w\in\Gamma$ such that $w(a_nb_na_nb_n\ldots
a_nb_nt')=a_nb_na_nb_n\ldots a_nb_na_n$. Then by
Lemma~\ref{lem_actions} there is $w'\in\Gamman$ such that
$w'(a_nb_na_nb_n\ldots a_nb_nt')=w(a_nb_na_nb_n\ldots
a_nb_nt')=a_nb_na_nb_n\ldots a_nb_na_n$. This proves transitivity of
$\Gamman$ on the levels of $\hat T^{(n)}$.
\end{proof}

%Note that if $n$ is even the reduction to the base case can be
%simplified. In this case for $\alpha'=\alpha_n^{n-1}$ and
%$\beta'=\beta_n^{n-1}$ we have
%\begin{equation}
%\begin{array}{lrl}
%\alpha'=&(a_n,b_n)(c_n,q_{n,n-2})\cdots(q_{n,n-3},q_{n,2n-4})&(\alpha',\alpha',\beta',\ldots,\beta'),\\
%\beta'=&(c_n,q_{n,n-2})\cdots(q_{n,n-3},q_{n,2n-4})&(\beta',\beta',\alpha',\ldots,\alpha').
%%\alpha_n^{n-1}=&(a_n,b_n)(c_n,q_{n,n-2})(d_n,q_{n,n-1})\cdots(q_{n,n-3},q_{n,2n-4})&(\alpha_n^{n-1},\alpha_n^{n-1},\beta_n^{n-1},\ldots,\beta_n^{n-1}),\\
%%\beta_n^{n-1}=&(c_n,q_{n,n-2})(d_n,q_{n,n-1})\cdots(q_{n,n-3},q_{n,2n-4})&(\beta_n^{n-1},\beta_n^{n-1},\alpha_n^{n-1},\ldots,\alpha_n^{n-1}),
%\end{array}
%\end{equation}
%Therefore the action of
%$\langle\alpha',\beta',\overline{(b_n,c_n)}\rangle<\Gamman$ on the
%tree consisting of the words over the alphabet
%$\{a_n,b_n,c_n,q_{n,n-2}\}$ is isomorphic to the action of $\Gamma$
%on $T$. This can be used analogously to the argument above to prove
%transitivity of $\Gamman$ on the levels of $\hat T^{(n)}$.

Finally, Theorem~\ref{thm_general_case} is derived from
Lemma~\ref{lem_gamman_trans} exactly in the same way as
Theorem~\ref{thm_main} is obtained from Lemma~\ref{lem_transit}.

\bibliographystyle{alpha}

\newcommand{\etalchar}[1]{$^{#1}$}
\def\cprime{$'$} \def\cprime{$'$} \def\cprime{$'$} \def\cprime{$'$}
  \def\cprime{$'$} \def\cprime{$'$} \def\cprime{$'$}

\noindent\texttt{savchuk@math.tamu.edu}

\noindent\texttt{yvorobet@math.tamu.edu}

\noindent
Department of Mathematics \\Texas A\&M University\\
College Station, TX 77843-3368 \\USA

\end{document}